\documentclass[a4paper,11pt,leqno]{amsart}
 
\usepackage{amssymb,hyperref}
\theoremstyle{plain}

\newtheorem{theo}{Theorem}[section]

\newtheorem{proposition}[theo]{Proposition}
\newtheorem{corollary}[theo]{Corollary}

\newtheorem{conjecture}[theo]{Conjecture}
\newtheorem{theorem}[theo]{Theorem}

\newtheorem{remark}[theo]{Remark}

\theoremstyle{definition}
\newtheorem{definition}[theo]{Definition}

\theoremstyle{remark}

 \def\RR{{\mathbb R}}
\def\CC{{\mathbb C}}
\def\ZZ{{\mathbb Z}}

\begin{document}

\title[Meromorphic geometric structures]{Meromorphic  almost  rigid geometric structures}

\author[S. Dumitrescu]{Sorin DUMITRESCU}

\address{D\'epartement de Math\'ematiques d'Orsay, 
\'equipe de Topologie et Dynamique,
Bat. 425, U.M.R.   8628  C.N.R.S.,
Univ. Paris-Sud (11),
91405 Orsay Cedex, France}
\email{Sorin.Dumitrescu@math.u-psud.fr}

\thanks{This work was  partially supported by  the ANR Grant Symplexe BLAN 06-3-137237}
\date{\today}

\setcounter{tocdepth}{1}

\maketitle


\begin{abstract} We study the local Killing Lie algebra of  meromorphic almost rigid geometric structures on complex manifolds. This leads to  classification results
for compact complex manifolds bearing holomorphic rigid geometric structures. \\ 
\end{abstract}

\section{Introduction}

Zimmer and Gromov have conjectured that ``big'' actions of Lie groups that preserve unimodular rigid geometric structures (for example, a pseudo-Riemannian metric or
an affine connection together with a volume form) are ``essentially classifiable''~\cite{DG, Gro,Zimmer1, Zimmer2, Zimmer3, Zimmer4}.

We try to adress here this general question in the framework of the complex geometry. We study holomorphic (and meromorphic) unimodular geometric structures on compact
complex manifolds. We show that in many cases  the holomorphic rigidity implies that the situation is classifiable (even without assuming any isometric group action). This is a survey
paper, but some results are new (see theorem~\ref{result}, corollary~\ref{corollary}  and theorem~\ref{meromorphic connections}).

In the sequel all complex manifolds are supposed to be smooth and connected. 

Consider  a complex $n$-manifold $M$ and, for all integers $r \geq 1$,   consider the associated bundle $R^r(M)$ of $r$-frames, which is a $D^r(\CC^n)$-principal bundle  over $M$, with $D^r(\CC^n)$-the (algebraic) group of $r$-jets at the origin of local biholomorphisms of  $\CC^n$ fixing $0$~\cite{DG,Gro, Benoist2, CQ, Feres, Zeghib}.

Let us consider, as in~\cite{DG,Gro}, the following 

\begin{definition} A {\it meromorphic geometric structure} (of order $r$)  $\phi$ on a complex manifold $M$ is a $D^r(\CC^n)$-equivariant map from  $R^r(M)$ to a quasiprojectif
variety $Z$ endowed with an algebraic action of $D^r(\CC^n)$. If $Z$ is an affine variety, we say that $\phi$ is of {\it affine type}.
\end{definition}

                    In holomorphic local coordinates $U$ on $M$, the geometric structure $\phi$ is given by a meromorphic function $\phi_{U} : U \to Z$.

                            Denote by $P \subset M$ the  analytic subset of poles of $\phi$ and by $M^*=M \setminus P$ the open dense subset of $M$ where $\phi$ is holomorphic. If $P$
                            is void,  then $\phi$ is a {\it holomorphic geometric structure} on $M$.
                            
   \begin{definition}  The  meromorphic geometric structure $\phi$ is called {\it almost  rigid} if, away from  an  analytic proper subset of $M$, $\phi$ is (holomorphic and) rigid
   in Gromov's sense~\cite{DG, Gro}. If on this open dense set $\phi$ uniquely determines a holomorphic  volume form, then $\phi$ is called {\it unimodular}.
   
                                     A holomorphic geometric structure on $M$  is called {\it rigid} if it is rigid in Gromov's sense on full $M$.
   \end{definition}
    
    Notice that our definiton of almost rigidity is not exactly the same as in~\cite{BF}.

  {\bf  Examples.} On a complex connected $n$-manifold $M$ the following meromorphic geometric structures are almost  rigid.
   
   \begin{itemize}
   
   \item  A meromorphic map $\phi$ from $M$ to a projective space $P^N(\CC)$, which is an embedding on some non trivial open subset of $M$. In this case $M$ is of algebraic  dimension
   $n$ (see section~\ref{meromorphic functions}). If $\phi$ is a holomorphic embedding, then $\phi$ is rigid.
   
   \item  A family $X_{1}, X_{2}, \ldots, X_{n}$ of meromorphic vector fields such that there exists an open (dense) subset of $M$ where the $X_{i}$ span the holomorphic tangent bundle $TM$.
   Any homogeneous manifold admits a holomorphic family of  such vector fields. They  pull-back  after  the blow up of a point (or a submanifold) in the homogeneous manifold
   in a meromorphic family of such vector fields.
   
  In the special case where  $X_{1}, X_{2}, \ldots, X_{n}$ are holomorphic and span $TM$ on full $M$, the corresponding geometric structure is called a {\it holomorphic parallelisation
  of the tangent bundle} and it is a (unimodular) holomorphic rigid geometric structure of affine type.
   
   \item  A meromorphic section $g$  of the bundle $S^2(T^{*}M)$ of complex quadratic forms such that $g$ is non degenerate on an open (dense) subset. If $g$ is holomorphic and non degenerate on $M$ then $g$ is a holomorphic rigid geometric structure of affine type called {\it  holomorphic Riemannian metric}. Up to a double cover of $M$, a holomorphic
   Riemannian metric is unimodular.
   
   A   holomorphic Riemannian metric has nothing to do with the more usual Hermitian metrics. It is in fact nothing but the complex version of (peudo-)Riemannian metrics. 
   Observe that  since complex quadratic forms  have no signature, there is here no distinction   between the Riemannian and pseudo-Riemannian cases. This  observation was at the origin of the nice use by F. Gau{\ss}  of the complexification technic of (analytic)  Riemannian metrics on surfaces, in order to find  conformal coordinates for them. Actually, the complexification of analytic Riemannian metrics leading to holomorphic ones, is becoming a standard trick (see for instance~\cite{Fra}).

   In the blow up process, the pull-back of a holomorphic Riemannian metric will stay non degenerate away from the exceptional set and will vanish on the exceptional set. 
   
   More general, if $g$ is a meromorpic section of $S^2(T^{*}M)$ non degenerate on an open dense set, then any complex manifold  $\tilde M$ bimeromorphic to   $M$ inherits   a
 similar   section $\tilde g$. Moreover, if $g$ is holomorphic, then $\tilde g$ is also holomorphic.   This comes from the fact that  indeterminacy points of a meromorphic map are of codimension
   at least 2~\cite{Ueno} (theorem 2.5) and from Levi's extension principle (see~\cite{Ueno} for the details).

   \item  Meromorphic affine or projective connections or meromorphic conformal structures in dimension $\geq 3$. For example, in local holomorphic coordinates $(z_{1}, \ldots, z_{n})$
   on the manifold $M$, a meromorphic affine connection $\nabla$ is determined by the meromorphic functions $\Gamma_{ij}^k$ such that $\nabla_{\frac{ \partial}{\partial z_{i}}}\frac{\partial}{\partial z_{j}}=\Gamma_{ij}^k \frac{\partial}{\partial z_{k}}$, for all $i,j,k \in \{1, \ldots, n\}$. 
   
   Affine connections are geometric structures of affine type, but projective connections or conformal structures are not of affine type.
   
   \item  If $\phi$ is a meromorphic  almost rigid geometric structure of order $r_{1}$ and $g$ is a meromorphic geometric structure of order $r_{2}$, then we can put together $\phi$ and $g$ in some  meromorphic  almost rigid geometric structure  of order $max(r_{1},r_{2})$, denoted by $(\phi,g)$~\cite{Gro, DG, Benoist2,CQ, Feres, Zeghib}.
   
   \item  The $s$-jet (prolongation) $\phi^{(s)}$ of a meromorphic almost rigid geometric structure $\phi$ is still a meromorphic almost rigid geometric structure~\cite{Gro, DG, Benoist2,CQ, Feres, Zeghib}.
    \end{itemize}
   
   Recall that  local biholomorphisms of $M$ which preserve a meromorphic geometric structure $\phi$ are called {\it local isometries}. Note that local isometries of $(\phi,g)$ are the local isometries of $\phi$ which preserve also $g$.
   
   The set of local isometries of a holomorphic rigid geometric structure $\phi$  is a Lie pseudogroup $Is^{loc}(\phi)$ generated by  Lie algebra of local vector fields called   {\it Killing Lie algebra}. If the Killing Lie algebra is  transitive on  $M$, then $\phi$ is called {\it locally homogeneous}.

   \section{Local isometries and meromorphic functions}    \label{meromorphic functions}
   
     The maximal number of $\CC$-linearly independent meromorphic functions on a complex manifold $M$ is called {\it the algebraic dimension} $a(M)$ of $M$. 
     
     Recall that a theorem of Siegel proves that a complex $n$-manifold $M$ admits at most $n$ linearly independent  meromorphic functions~\cite{Ueno}.
      Then $a(M) \in \{0, 1, \ldots, n\}$ and for algebraic manifolds $a(M)=n$.

   We will say that two points in  $M$ are in the same {\it fiber of the algebraic reduction} of $M$ if any meromorphic  function on $M$ takes the same value at  the two points. There
   exists some open dense set in  $M$ where the fibers of the algebraic  reduction are the fibers of a holomorphic fibration  on an algebraic manifold of dimension $a(M)$ and any 
   meromorphic function on $M$ is the pull-back of a meromorphic function on the basis~\cite{Ueno}.

   The following theorem shows  that the fibers of the algebraic reduction are in the same orbit of the pseudogroup of local isometries for any meromorphic almost  rigid
   geometric structure on $M$. This is a meromorphic  version of the celebrated Gromov's  open-dense orbit theorem~\cite{DG,Gro} (see also~\cite{Benoist2,CQ, Feres, Zeghib}).

\begin{theorem} \label{result}
Let $M$ be a  connected complex manifold of dimension $n$ which admits a meromorphic almost  rigid geometric structure $\phi$. 
Then,  there exists a nowhere dense analytic subset $S$ of $M$, such that $M \setminus S$ is $Is^{loc}(\phi)$-invariant and the orbits of $Is^{loc}(\phi)$ in $M \setminus S$
are the fibers of a holomorphic fibration of constant rank. The dimension of the fibers  is $\geq n-a(M)$, where $a(M)$ is the algebraic dimension
of $M$.
\end{theorem}

Recall that $g'=(\phi,g)$ is still  almost rigid for any meromorphic (not necessarily almost rigid)  geometric structure $g$ on $M$. This yields to the following:

\begin{corollary}~\label{corollary} If $g$ is a meromorphic  geometric structure on $M$,  away from a  nowhere dense  analytic subset,  the $Is^{loc}(g)$-orbits are of dimension
$\geq n-a(M)$.

In particular, if $a(M)=0$, then $g$ is locally homogeneous outside a nowhere dense  analytic subset of $M$ (the Killing Lie algebra $\mathcal G$ of $g'=(\phi,g)$ is transitive
on an open dense set). Moreover, if $\phi$ is unimodular and $\mathcal G$ is unimodular and simply transitive, then $g$ is locally homogeneous on the full open dense set where $\phi$ and $g$ are holomorphic
and $\phi$ is rigid. In this case, if  $\phi$ is unimodular  holomorphic and  rigid and $g$ is holomorphic, then $g$ is locally homogeneous on $M$.
\end{corollary}

\begin{proof}
 For each positif  integer  $s$  we consider the $s$-jet $\phi^{(s)}$ of  $\phi$. This is a $D^{(r+s)}(\CC^n)$-equivariant meromorphic map 
    $R^{(r+s)}(M) \to Z^{(s)},$ where $Z^{(s)}$ is the algebraic variety of the $s$-jets at the origin of  holomorphic maps from $\CC^n$ to $Z$. One can find the expression of the (algebraic)
    $D^{(r+s)}(\CC^n)$-action on $Z^{(s)}$ in~\cite{CQ, Feres}.
    
    Since $\phi$ is almost rigid, there exists a nowhere dense  analytic subset $S$ in $M$ and a positive integer  $s$ such that  two points $m,m'$ in  $M \setminus S$ are in the same orbit 
   of $Is^{loc}(\phi)$ if and only if $\phi^{(s)}$ sends the  fibers  of $R^{(r+s)}(M)$ above $m$ and $m'$
    on the same  $D^{(r+s)}(\CC^n)$-orbit in  $Z^{(s)}$~\cite{Gro, DG}.
   
     Rosenlicht's  theorem (see~\cite{Po}) shows that there exists  a $D^{(r+s)}(\CC^n)$-invariant stratification 
     
     $$Z^{(s)}=Z_0 \supset \ldots \supset Z_l,$$
      such that 
     $Z_{i+1}$ is Zariski closed in  $Z_i$, the quotient of $Z_{i} \setminus Z_{i+1}$ by $D^{(r+s)}(\CC^n)$ is a complex manifold  and 
     rational $D^{(r+s)}(\CC^n)$-invariant functions  on  $Z_i$ separate
      orbits in  $Z_i \setminus Z_{i+1}$.
     
     Take the open dense $Is^{loc}(\phi)$-invariant subset $U$ of $M$, where $\phi^{(s)}$ is of constant rank and the 
      image of  $R^{(r+s)}(M)|_{U}$ through  $\phi^{(s)}$ is contained in the  maximal subset  $Z_i \setminus Z_{i+1}$ of the stratification which intersects the image of $\phi^{(s)}$. Then the orbits of $Is^{loc}(\phi)$ in $U$
      are the fibers of a fibration of constant rank.
      
      If $m$ and $n$ are two points in $U$ which are not in the same $Is^{loc}(\phi)$-orbit, then the corresponding  fibers of $R^{(r+s)}(M)|_{U}$ are sent by $\phi^{(s)}$  on two
      distinct $D^{(r+s)}(\CC^n)$-orbits in $Z_i \setminus Z_{i+1}$. By Rosenlicht's theorem there exists a $D^{(r+s)}(\CC^n)$-invariant  rational function $F : Z_i \setminus Z_{i+1} \to \CC,$
      which takes distincts values at these two orbits.

     The meromorphic function  $F \circ \phi^{(s)} : R^{(r+s)}(M) \to \CC$ is $D^{(r+s)}(\CC^n)$-invariant and descends in a $Is^{loc}(\phi)$-invariant meromorphic function on 
     $M$ which takes distincts values at  $m$ and at $n$.
     
    Consequently, the complex codimension in $U$ of the $Is^{loc}(\phi)$-orbits is $\leq a(M)$, which finished the proof.
    \end{proof}

We deduce now  the corollary:

\begin{proof}

It is  convenient to put together $\phi$ and $g$ in some extra geometric structure $g'=(\phi, g)$. Now $g'$ is a meromorphic almost
rigid geometric structure and  theorem~\ref{result} shows  that the complex dimension of the orbits of $Is^{loc}(g')$ is $\geq n-a(M)$. As each  $Is^{loc}(g')$-orbit is  contained in a $Is^{loc}(g)$-orbit, the result follows also for $g$.

If $a(M)=0$, the Killing Lie algebra $\mathcal G$ of $g'$ is transitive on an maximal open dense subset $U$ of $M$. Suppose now $\phi$ is unimodular and also $\mathcal G$ is unimodular. The other inclusion being trivial, it is enough  to show that $U$ is the maximal subset of  $M$ where $\phi$ and $g$ are holomorphic and $\phi$ is rigid. 

Take a point $m$ in the previous subset. We want to show that $m$ is in $U$. Since  $g'=(\phi,g)$ is holomorphic and rigid in the neighborhood of $m$, it follows that $m$
admits an open neighborhood $U_{m}$ in $M$ such that any local holomorphic Killing field of $g'$ defined on a connected open subset in $U_{m}$ extends on $U_{m}$~\cite{Nomizu,Amores}.

Since $\mathcal G$ acts transitively on $U$, choose local linearly independent Killing fields  $X_{1}, \ldots, X_{n}$ on a  connected open set included  in $U \cap U_{m}$. As
$\phi$ is unimodular, it determines a holomorphic volume form on $U_{m}$ (if necessary we restrict to  a smaller $U_{m}$) which is preserved by $Is^{loc}(g')$. But $Is^{loc}(g')$
acts transitively on $U \cap U_{m}$ and $\mathcal G$ is supposed to be unimodular. This implies that the function $vol(X_{1}, \ldots, X_{n})$ is $\mathcal G$-invariant and, consequently, a (non vanishing) constant on $U \cap U_{m}$.

On the other hand $X_{1}, \ldots, X_{n}$ extend in some holomorphic Killing fields $\tilde X_{1}, \ldots, \tilde X_{n}$ defined on full $U_{m}$. The holomorphic function
$vol(\tilde X_{1},  \ldots, \tilde X_{n})$ is a non vanishing constant on $U_{m}$: in particular, $\tilde X_{1}(m), \ldots, \tilde X_{n}(m)$ are linearly independent. We proved that
$\mathcal G$ acts transitively in the neighborhood of $m$ and thus $m \in U$.
\end{proof}

\section{Classification results}

\subsection*{Parallelisations of the tangent bundle}

We begin with the classification of compact complex manifolds which admits a holomorphic parallelisation of the tangent bundle~\cite{Wang}.

\begin{theorem}(Wang) Let $M$ be a compact connected complex manifold with a holomorphic parallelisation of the tangent bundle. Then $M$ is a quotient $\Gamma \backslash G$, where
$G$ is a simply connected complex Lie group and $\Gamma$ is a uniform lattice in it.

Moreover, $M$ is Kaehler if and only if $G$ is abelian (and $M$ is a complex torus).
\end{theorem}

\begin{proof}

Take $X_{1}, X_{2}, \ldots, X_{n}$ global holomorphic vector fields on $M$ which span $TM$. Then, for all $1\leq i,j \leq n$, we have 

$$\lbrack X_{i}, X_{j} \rbrack =f_{1}^{ij} X_{1}  + f_{2}^{ij} X_{2} +\ldots+ f_{n}^{ij}X_{n},$$ 

with $f_{k}^{ij}$ holomorphic functions on $M$. Since $M$ is compact, these functions has to be constant by the maximum principle and, consequently, 
$X_{1}, X_{2},  \ldots, X_{n}$ generate a n-dimensional Lie algebra  $\mathcal G$ which acts simply transitively on $M$. By Lie's theorem there exists a unique connected simply connected complex Lie group $G$ corresponding to $\mathcal G$.

In particular, the holomorphic parallelisation is locally homogenous, locally modelled on the parallelisation given by translations-invariant vector fields on the Lie group $G$.

Since $M$ is compact, the $X_{i}$ are complete and they define a holomorphic simply transitive action of $G$ on $M$. Hence $M$ is a quotient of $G$ by a cocompact discrete
subgroup of $G$.

Assume now $M$ is Kaehler. Then, any holomorphic form on $M$ has to be closed~\cite{GH}. Consider $\omega_{1}, \omega_{2}, \ldots, \omega_{n}$ the dual basis in respect to
$X_{i}$. The one forms $\omega_{i}$ are holomorphic and translations-invariants on $G$. The Lie-Cartan formula

$$d\omega_{i}(X_{j},X_{k})=-\omega_{i}(\lbrack X_{j},X_{k} \rbrack),$$

shows that the one forms $\omega_{i}$ are all closed if and only if $G$ is abelian and thus $M$ is a complex torus.
\end{proof}

\subsection*{Holomorphic Riemannian metrics}

As in the real case,  a  holomorphic Riemannian metric on $M$ gives rise to  a covariant differential calculus, i.e.    a  Levi-Civita (holomorphic) affine 
connection, and to geometric features:   curvature tensors, geodesic (complex) curves~\cite{Le, Leb}.

Locally,  a holomorphic Riemannian metric has the form 
$\Sigma g_{ij}(z)dz_i dz_j$, where $(g_{ij}(z))$ is a complex inversible   symmetric matrix depending holomorphically on $z$.
The standard   example is that of  the global flat  holomorphic Riemannian metric  $dz_{1}^2+dz_{2}^2+ \ldots +dz_{n}^2$ on $\CC^n$. This metric is translations-invariant and thus 
descends  to any quotient of $\CC^n$ by a lattice. Hence complex tori  possess  (flat) holomorphic Riemannian metrics.  This is however a very  special situation since, contrary to 
 real case, {\it only  few    compact complex manifolds admit holomorphic Riemannian metrics.} In fact,  Yau's proof of the Calabi conjecture shows that, up to finite unramified  covers,  complex tori are the only compact {\it Kaehler}  manifolds admitting holomorphic Riemannian metrics~\cite{IKO}. 
 
 However, very interesting examples, constructed by Ghys in~\cite{Ghys}, do exist on $3$-dimensional complex non Kaehler manifolds  and deserve classification. Notice that   parallelisable
 manifolds, bear holomorphic Riemannian metrics coming from left invariant 
 holomorphic Riemannian metrics on $G$ (which can be constructed  by left translating  any complex non-degenerate quadratic form defined on the Lie algebra ${\mathcal G}$.)

 Ghys's examples of $3$-dimensional compact complex manifolds endowed with holomorphic Riemannian metrics are obtained by deformation of the complex structure
 on parallelisable manifolds $\Gamma \backslash  SL(2, \CC)$~\cite{Ghys}. They are {\it non standard, meaning they do not admit  parallelisable manifolds as finite unramified covers}. Those non standard examples will be described further on.
 
 A first obstruction to the existence of a holomorphic Riemannian metric on a compact complex manifold is  its  first Chern class. Indeed, a holomorphic Riemannian
metric on $M$ provides an isomorphism between $TM$ and $T^{*}M$. In particular, the canonical bundle $K$ is isomorphic to the anti-canonical bundle $K^{-1}$
and $K^2$ is trivial. This means that the first Chern class of $M$ vanishes and, up to a double unramified cover, $M$ possesses a holomorphic volume form.

 The following proposition describes  holomorphic Riemannian metrics on  parallelisable manifolds:
 
 \begin{proposition}  Let $M=\Gamma \backslash G$ a compact parallelisable manifold, with $G$ a simply connected complex Lie group and $\Gamma$ a uniform lattice in $G$.
 Then,  any holomorphic Riemannian metric $g$ on $M$ comes from a non degenerate  complex quadratic form on the Lie algebra $\mathcal G$ of $G$. In particular, the pull-back of $g$ is right invariant on the universal cover $G$ (and  $g$ is locally homogeneous on $M$).
 
 Moreover, any compact parallelisable $3$-manifold admits a holomorphic Riemannian metric of constant sectional curvature. The metric is flat exactly when $G$ is solvable.
 \end{proposition}
 
 \begin{proof}
 
 Consider $X_{1}, X_{2}, \ldots, X_{n}$ the fundamental vector fields corresponding to the simply transitive $G$-action on $M$. Let $g$ be a holomorphic Riemannian metric on $M$ and denote also by $g$ the associated complex
symmetric  bilinear form. Then $g(X_{i},X_{j})$ is a holomorphic function on $M$ and thus constant, for all $1 \leq i,j \leq n.$ This implies that $g$ comes from a right-invariant holomorphic Riemannian metric on $G$.
 
 Assume now  $G$ is a simply connected connected complex unimodular Lie group of dimension $3$. We have only four such Lie groups: $\CC^3$,
 the complex Heisenberg group, the complex SOL group and $SL(2, \CC)$~\cite{Kir}. Note that the group $SOL$ is the complexification of  the  affine isometry group of the Minkowski plane $\RR^{1, 1}$ or equivalently the isometry group of $\CC^2$ endowed with its flat holomorphic Riemannian  metric.
 
We begin with the case $G=\CC^3$.  We have seen that $\CC^3$ admits a  flat translations-invariant holomorphic Riemannian metric. The isometry group of this metric is $O(3, \CC)  \ltimes \CC^3$. We recall here
 that $O(3, \CC)$ and $SL(2, \CC)$ are locally isomorphic.
 
 Consider now the case $G=SL(2, \CC)$. The Killing form of the Lie algebra $sl(2, \CC)$ is a  non degenerate complex quadratic form which endows $SL(2, \CC)$ with a left invariant holomorphic Riemannian 
 metric of constant sectional curvature. Since the Killing quadratic form is invariant by the adjoint representation the isometry group contains also all the right translations.
 In fact the connected component of the isometry group is $SL(2, \CC) \times SL(2, \CC)$ acting by left and right translations.
 
 It is an easy exercice to exhibit in the isometry group $O(3, \CC) \ltimes \CC^3$ of the flat holomorphic Riemannian space, copies of the complex Heisenberg group and of the complex SOL group
 which act simply and transitively. Thus  the  flat holomorphic Riemannian space  also admits models which are given by left invariant metrics on the Heisenberg group and on
 the SOL group.  One can get  explicit expression  of these holomorphic Riemannian metrics by complexification of flat left invariant Lorentz metrics on the real Heisenberg and SOL
 groups~\cite{Ra, Rah}.
 \end{proof}
 
 {\bf Ghys's non standard examples.}  As above, for  any co-compact lattice $\Gamma$ in $SL(2, \CC)$, the quotient 
  $M=  \Gamma \backslash SL(2, \CC)$  admits a holomorphic Riemannian metric  of non-zero constant sectional curvature.
  It is convenient to consider   $M$ as  a quotient of $S_3=O(4, \CC) /  O(3, \CC) =SL(2, \CC) \times SL(2, \CC) /  SL(2, \CC)$ by 
  $\Gamma$,  seen as a subgroup of $SL(2, \CC) \times SL(2, \CC)$ by the 
  trivial embedding $\gamma \in \Gamma \mapsto ( \gamma,1) \in SL(2, \CC) \times SL(2, \CC)$.

New interesting examples of manifolds admitting holomorphic Riemannian metrics of non-zero constant sectional curvature have been constructed in \cite{Ghys}  by deformation of this embedding of $\Gamma$. 

Those  deformations  are constructed   choosing   a morphism 
 $u: \Gamma \to SL(2, \CC)$ and considering the embedding 
 $ \gamma \mapsto ( \gamma, u(\gamma))  $.
 Algebraically, the  action is given 
 by:
 $$(\gamma,m) \in \Gamma \times SL(2,\CC)  \to \gamma m u(\gamma^{-1}) \in SL(2,\CC).$$

 It is proved in \cite{Ghys} that, for $u$ close enough to the trivial morphism,  $\Gamma$ acts properly (and freely)  on $S_3 (\cong SL(2, \CC))$ such that the quotient $M_u$
  is  a complex compact manifold (covered by $SL(2, \CC)$) admitting a holomorphic Riemannian metric of non-zero constant sectional curvature.  {\it In general, these examples do
  not admit parallelisable manifolds as finite covers}.

   Note that  left-invariant holomorphic
  Riemannian metrics on $SL(2, \CC)$ which are not right-invariant, in general, will not  descend on $M_{u}$.

Let us notice   that despite this systematic study in 
\cite{Ghys}, there are  still many open questions regarding these examples (including the question of completeness).  A real version of this study is in~\cite{Kulkarni-Raymond, Goldman, Salein}.

{\bf Dimension 3}. The classification of complex compact manifolds admitting holomorphic Riemannian metrics is simple in complex dimension $2$ (see section~\ref{surfaces}). 
 An important step
 toward the classification in  dimension $3$ was  made in~\cite{Dum1} with the following result:

\begin{theorem}  \label{locally homogeneous}  
Any holomorphic Riemannian metric on a compact connected complex $3$-manifold is locally homogeneous. More generally, 
 if  a compact connected complex $3$-manifold $M$ admits a   holomorphic Riemannian metric, then any  holomorphic geometric structure of affine type on $M$
is locally homogeneous. 
\end{theorem}

Observe that the previous result is trivial in  dimension $2$, since the sectional curvature is a holomorphic function and thus constant on compact complex surfaces: this
implies the local homogeneity~\cite{Wo}. In dimension $3$, the sectional curvature will be, in general,  a non constant meromorphic function on the $2$-grassmanian of the holomorphic tangent space with poles on the degenerate planes.

Thanks to theorem~\ref{locally homogeneous}, our manifold $M$ is locally modelled on a $(G, G/I)$-geometry  in Thurston's sense~\cite{Th}, where $I$ is a closed subgroup of the Lie group $G$ such that the $G$-action on  $G/I$ preserves some holomorphic Riemannian metric (notice that the local {\it Killing Lie algebra} of the holomorphic Riemannian metric is the Lie algebra of $G$). In this contexte we have a {\it developping map} from the universal
cover of $M$ into $G/I$ which is a local diffeomorphism and which is equivariant in respect to the action of the fundamental group on $\tilde M$ by deck transformations and on $G/I$ via the {\it holonomy morphism} $\rho: \pi_1(M) \to G$~\cite{Th}.

Recall that the $(G,G/I)$-geometry is called {\it complete} if the developping map is a diffeomorphism and, consequently, $\Gamma= \rho (\pi_1(M))$ acts properly on $G/I$  such that $M$ is a compact quotient  $\Gamma  \backslash G /I$. 

A second step was achieved in~\cite{DZ} were, in a commun work with Zeghib,  we proved the following result which can be seen, in particular, as a completeness result
in the case where $G$ is solvable.

\begin{theorem} \label{with Zeghib}
Let $M$ be a compact connected complex $3$-manifold which admits a (locally homogeneous) holomorphic Riemannian metric $g$. Then: 

(i) If the Killing Lie algebra of $g$ has a non trivial  semi-simple part, then it  preserves some  holomorphic Riemannian metric on $M$ with  constant sectional curvature.

(ii) If the Killing Lie algebra of $g$ is solvable, then, 
up to  a finite  unramified cover,
 $M$   is a quotient either of  the complex Heisenberg group,  or of the complex $SOL$ group by a lattice.
\end{theorem}

\begin{remark} If $g$ is flat, its Killing Lie algebra corresponds to $O(3, \CC) \ltimes \CC^3$, which has non trivial semi-simple part. Thus, flat holomorphic Riemannian metrics on complex tori are  part  of  point (i) in  the main  theorem.
\end{remark}

The point (ii) of the previous theorem is not only about  completeness, but also gives a  rigidity result in Bieberbach's sense~\cite{Wo}: $G$ contains  a 3-dimensional closed subgroup  $H$ (either isomorphic to  the complex
Heisenberg group, or to the complex $SOL$ group) which acts simply and transitively (and so identifies) with $G/I$ and (up to a finite index)  the image $\Gamma$ of the holonomy morphism  lies in $H$. It follows that, up to a finite cover,  $M$ identifies with $\Gamma \backslash H$.

Since the Heisenberg and SOL groups admit left  invariant flat holomorphic Riemannian metrics this leads  to:

\begin{corollary}  \label{principal}
 If a compact connected complex  3-manifold  $M$ admits a holomorphic Riemannian metric, then, up to  a finite   unramified cover,  $M$   admits a holomorphic Riemannian
 metric of constant sectional curvature. 
\end{corollary}

Theorem~\ref{with Zeghib}  does not end the story, even in dimension $3$, essentially because of remaining {\it completeness}  questions, and those on the algebraic structure of the fundamental group. 

{\it  It remains to classify the compact complex $3$-manifolds endowed with a holomorphic Riemannian metric of  constant sectional curvature}.  

{\bf Flat case.} In this case $M$ admits a  $(O(3, {\mathbb C}) \ltimes {\mathbb C}^3, {\mathbb C}^3)$-geometry.   The challenge  remains:

1) {\it Markus conjecture:}  Is $M$ complete?

2) {\it Auslander conjecture: }  Assuming $M$ as above, is $\Gamma$ solvable? 

Note that these questions are settled in the setting of (real) flat Lorentz manifolds~\cite{Car, Fried-Goldman}, but unsolved  for general (real) pseudo-Riemannian metrics. The real part
of the holomorphic Riemannian metric is a (real) pseudo-Riemannian metric of signature $(3,3)$ for which both previous conjectures are still open.

{\bf Non flat case}. In this case $G=SL(2, \CC) \times SL(2, \CC)$ and $I=SL(2, \CC)$ is diagonally embedded in the product. The completeness of this
geometry on compact complex manifolds is still an open problem, despite a local result of Ghys~\cite{Ghys}. Recall that the real analogous of this problem, i.e. the completeness
of compact manifolds endowed with Lorentz metrics of negative constant sectional curvature, was solved in~\cite{Klingler2}, but the proof cannot generalize  to other signatures.

{\bf Higher dimension.} One interesting 
 problem in differential geometry is to decide if a given homogeneous space $G/I$ possesses or not a compact quotient. A more general related question is to decide
 if there exist compact manifolds locally modelled on $(G,G/I)$ (see, for instance,~\cite{Benoist,  BL,  Kobayashi}).

 The case $I = 1$, or more generally $I$ compact,  reduces to the classical question of existence 
 of co-compact lattices in Lie groups.  For  homogeneous  spaces  of non-Riemannian  type (i.e. $I$ non-compact) the problem is much harder. 
 
The  case  
 $S_n = O(n+1, \CC)/ O(n, \CC)$ is a geometric situation where these questions can be tested. It turns out that compact quotients of  $S_n$ are known to exist only for $n=1,3$ or $7$. We discussed the case $n= 3$ above, and the existence of a compact quotient of $S_7$
 was proved  in 
 \cite{Kobayashi}.  Here, we dare ask with \cite{Kobayashi}:
 
 \begin{conjecture}\cite{Kobayashi} $S_n$ has no  compact quotients,  for $n \neq 1, 3, 7$.

 \end{conjecture}

 A stronger version of this question was  proved in \cite{Benoist}  for $S_n$, if $n$ has the form 
  $4m+1$.
  
  Keeping in mind our geometric approach, we generalize the question 
  to manifolds locally modelled on $S_n$. More exactly:

  \begin{conjecture}\cite{DZ} A compact complex manifold endowed with a holomorphic Riemannian metric
  of constant non-vanishing curvature is complete. 
   In particular,  such a manifold has dimension 3 or 7. 
  \end{conjecture}

\section{Applications to  simply connected manifolds}

Remark first that theorem~\ref{locally homogeneous} has the following direct consequence:

\begin{corollary} A compact connected simply connected $3$-manifold doesn't bear any holomorphic Riemannian metric.
\end{corollary}

\begin{proof} Assume by contradiction $(M,g)$ as in the statement of the corollary. Then theorems~\ref{locally homogeneous} implies that $g$ is locally homogeneous.  Since  $M$ is simply connected, the local Killing fields of the Killing algebra $\mathcal G$ extend on full $M$: the unique  connected simply connected complex Lie group $G$ associated to $\mathcal G$  acts isometrically and transitively on $M$. Then $M$ is a homogeneous space $G/H$. Moreover, up to a double cover,   $G/H$ admits a holomorphic volume form $vol$ coming from the holomorphic Riemannian metric.

Take $X_{1}, X_{2}, X_{3}$ three global Killing fields on $M$ which are linearly independent at some point. Since $vol(X_{1}, X_{2}, X_{3})$ is a holomorphic function on $M$, it
is a non zero  constant and, consequently, $X_{1}, X_{2}, X_{3}$ are linearly independent on $M$. Hence Wang's theorem  implies  that $M$ is a quotient of a three dimensional
connected simply connected complex Lie group $G_{1}$ by a discrete subgroup. Since $M$ is simply connected, this discrete subgroup has to be trivial and $M$ identifies with
$G_{1}$. But there is no compact  simply connected complex Lie group: a contradiction.
\end{proof}

\begin{theorem} \label{simply connected} A compact connected  simply connected complex $n$-manifold without  nonconstant meromorphic functions doesn't bear any holomorphic unimodular rigid geometric
structure.
\end{theorem}

\begin{proof}

Assume, by contradiction, that  $(M, \phi)$ verifies   the hypothesis. Since $a(M)=0$, theorem~\ref{result} implies  $\phi$ is locally homogeneous on an open dense set $U$. As $M$
is simply connected, elements  in  the Killing algebra $\mathcal G$ extend on full $M$: the unique  connected simply connected complex Lie group $G$ associated to $\mathcal G$  acts isometrically on $M$ with an open dense orbit. The open dense orbit $U$  identifies with a homogeneous space $G/H$, where $H$ is a closed subgroup of $G$.

Consider  $X_{1}, X_{2}, \ldots, X_{n}$  global Killing fields on $M$ which are linearly independent at some point of the open orbit $U$. As before, $vol(X_{1},X_{2}, \ldots, X_{n})$
is a non zero constant, where  $vol$ is the holomorphic volume form associated to $\phi$. Thus the $X_{i}$ give a holomorphic parallelisation of $TM$ and Wang's theorem enables
us to conclude as in the previous proof.
\end{proof}

For non unimodular rigid geometric structures we have the following less precise:

\begin{theorem} \label{dim3} \cite{D1} Let $M$ be a compact connected simply connected complex $n$-manifold without non constant meromorphic functions and admitting a holomorphic rigid geometric structure $\phi$. Then $M$ is an equivariant compactification of $\Gamma \backslash G$, where $\Gamma$ is a discrete non cocompact subgroup in  a complex Lie group $G$.
\end{theorem}

\begin{proof}
Since $a(M)=0$, theorem~\ref{result} implies  $\phi$ is locally homogeneous on an open dense set $U$. As before, the extension property of local Killing fields implies $U$ is 
a complex homogeneous space $G/H$, where $G$ is a connected simply connected complex Lie group and $H$ is a closed subgroup  in $G$.

We show now that $H$ is a discrete subgroup of $G$. Assume by contradiction the Lie algebra of $H$  is not trivial.
Take at any point $u \in U$, the isotropy subalgebra $\mathcal{H}_{u}$ (i.e. the Lie subalgebra of  Killing fields vanishing at $u$).
Remark that $\mathcal{H}_{gu}=Ad(g) \mathcal{H}_{u}$, for any $g \in G$ and $u \in U$, where $Ad$ is the adjoint representation.

The map $u \to \mathcal{H}_{u}$ is a meromorphic map from $M$ to the grassmanian of $d$-dimensional vector spaces in $\mathcal G$. But $M$ doesn't admit any non trivial
meromorphic function and this map has to be constant. It follows that $\mathcal {H}_{u}$ is $Ad(G)$-invariant and $H$ is a normal subgroup of $G$: a contradiction, since the $G$-action
on $M$ is faithfull. Thus $G$ is of dimension $n$ and $H$ identifies to a lattice $\Gamma$ in $G$.

As $M$ is simply connected,  $U$ has to be strictly included in $M$ and $M$ is an equivariant compactification of $\Gamma  \backslash G$.
\end{proof}

We don't know if such equivariant compactifications of $ \Gamma \backslash G$ admit equivariant holomorphic rigid geometric structures, but the previous result has the following application.
For the sake,  recall that it is still an open question to know if the $6$-dimensional real sphere $S^6$ bears complex structures. In this context, we have the following:

\begin{corollary} If $S^6$ admits a complex structure $M$, then $M$ doesn't bear any holomorphic rigid geometric structure.
\end{corollary}

\begin{proof} The starting point of the proof is a result of~\cite{CDP} where it is proved that $M$ doesn't admit any non constant meromorphic functions. If $M$
bears any holomorphic rigid geometric structure, then theorem~\ref{dim3}
implies that $M$ is an equivariant compactification of a homogeneous space. This is in contradiction with the main theorem  of~\cite{HKP}.
\end{proof}

As for {\it Kaehler manifolds} we have proved in~\cite{D1} the following more precise:

\begin{theorem}  \label{Kaehler}  Let $M$ be a compact connected Kaehler manifold endowed with a holomorphic unimodular rigid  geometric structure $\phi$ of affine type. Then, up to a finite unramified cover,
$M$ is a complex torus (quotient of $\CC^n$ by a lattice) and $\phi$ is translations-invariant.
\end{theorem}
 
 The proof is done in two steps. First we prove  that $\phi$ is  locally homogeneous. Then we use a splitting  theorem~\cite{Beauville}  which asserts  that such compact Kaehler manifolds with a holomorphic volume form  are biholomorphic to a direct product of a complex torus and a compact {\it simply connected} Kaehler manifold with a holomorphic volume form.  Starting with $\phi$ and using the product structure, we construct a holomorphic unimodular rigid geometric structure on the simply connected factor which is locally homogeneous. We conclude as in the proof of theorem~\ref{simply connected}
 that the simply connected factor is trivial.

In a recent work, McKay~\cite{McKay}  proved a similar result for holomorphic Cartan parabolic geometries (for example, conformal structures). One can also find a classification
of certain  holomorphic $G$-structures of order one on {\it uniruled} projective  manifolds in~\cite{HM}.

\section{Applications to complex surfaces}  \label{surfaces}

\begin{theorem}    Let $S$ be a complex compact surface endowed with a holomorphic unimodular rigid geometric structure $\phi$ of affine type. Then the Killing Lie algebra of $\phi$ is non trivial.
\end{theorem}

The assumption on $\phi$ to be affine is essential:

\begin{remark}  If $S$ is an algebraic compact complex surface with trivial canonical bundle (for example, a complex torus or an algebraic K3 surface),  the geometric structure given by a  holomorphic volume forme on $S$ together with  a holomorphic embedding of $S$ in a complex projective space doesn't admit any non trivial local isometry. However
this geometric structure is not of {\it affine type}.
\end{remark}

\begin{proof}
The proof  is a simply corollary of theorems~\ref{Kaehler} and~\ref{result}. Indeed, theorem~\ref{result} implies the Killing Lie algebra of $\phi$ is trivial only if the algebraic  dimension
$a(S)$ equals $2$. But in this case $S$ is algebraic~\cite{Barth-Peters} and  thus Kaehler. Then  theorem~\ref{Kaehler} applies and the Killing Lie algebra is transitive on $S$ and hence of dimension
at least $2$.
\end{proof}

\begin{theorem}
Let $S$ be a mimimal complex compact surface which is not biholomorphic to  a non algebraic $K3$ surface or to a non affine  Hopf surface. Then $S$ admits holomorphic  rigid geometric structure.
\end{theorem}

\begin{remark} By definition,  a complex algebraic manifold admits an embedding in a complex projectif space, which was seen to be a holomorphic rigid geometric structure
(of order zero).
\end{remark}

\begin{proof}

By the previous remark, it remains to consider the case of non algebraic complex surfaces. Then the   Enriques-Kodaira classification shows that, up to a finite unramified cover,
any minimal non algebraic complex compact surface is biholomorphic to one of the following complex surfaces: a complex tori, a Hopf surface, an Inoue surface, a K3 surface, an elliptic principal bundle over an elliptic curve or an elliptic principal bundle over a Riemann surface of genus $g \geq 2$, with odd first Betti number (see~\cite{Barth-Peters}, p. 244).

Or, it is known that  complex tori, Inoue surfaces, affine  Hopf surfaces,  elliptic principal bundles over  elliptic curves and elliptic principal bundles over a surface of genus $g \geq 2$ with odd first Betti number admit flat holomorphic affine connections~\cite{IKO, Klingler, Maehara, Suwa, Vitter, Wall}.
\end{proof}

The classification of complex surfaces admitting holomorphic affine connections (see~\cite{IKO, Klingler, Maehara, Vitter, Wall}) implies the following result:

\begin{theorem} Let $S$ be a complex compact surface admitting a holomorphic  unimodular affine connection. Then, up to a finite unramified cover, either $S$ is a complex torus and the connection is translations-invariant,  or $S$ is an elliptic principal bundle over an elliptic curve and the connection is locally modelled on a translations-invariant connection
on a complex torus.
\end{theorem}

Since any complex manifold endowed with a  holomorphic Riemannian metric inherits  a  holomorphic (unimodular) affine (Levi-Civita) connection~\cite{Le, Leb}, this easily implies the
following~\cite{Dum1}:

\begin{corollary} Let $S$ be a complex compact surface admitting a holomorphic Riemannian metric $g$. Then, up to a finite unramified cover, $S$ is a complex torus
and $g$ is   (flat) translations-invariant.
\end{corollary}

We give here a direct proof of the flatness of $g$, without using the Enriques-Kodaira classification.

\begin{proof} The constant sectional curvature is a holomorpic function on $S$ and hence  it is constant. One can multiply a holomorphic Riemannian metric by a complex constant $\lambda$ which induces a multiplication by $\lambda^{-2}$ of its sectional curvature. Therefore, only 
 the vanishing or not (but not the sign)  of the curvature
 is relevant. 
 
 Assume, by contradiction, that $g$ is of constant non zero curvature. Then, up to rescaling,  $g$ is locally modelled to the following model. Consider the adjoint representation
 of $SL(2, \CC)$ into its Lie algebra $sl(2, \CC)$. This $SL(2, \CC)$-action preserves the Killing quadratic form $q$, which is a non degenerate complex quadratic form.
 Choose $x \in sl(2, \CC)$ a vector of unitary $q$-norm and consider the orbit of $x$. This orbit is a homogeneous space
 $SL(2, \CC)/I$, where $I$ a one-parameter semi-simple subgroup of $SL(2, \CC)$, on which the restriction of the Killing form induces a two-dimensional  $SL(2, \CC)$-invariant  complete holomorphic Riemannian metric. This model has constant non zero sectional curvature and $g$ has to be locally modelled on it.
 
 We prove now that there is no {\it compact} complex surface locally modelled on $SL(2, \CC)/I$.
 
 Observe that $x$ induces on $SL(2,\CC)/I$ a $SL(2, \CC)$-invariant vector field $X$ and $g(X, \cdot)$ induces on $SL(2, \CC)/I$ an invariant one form $\omega$. Remark  that 
 $d \omega$ is a volume form. Indeed, $d\omega(Y,Z)=-g(X, \lbrack Y, Z \rbrack)$, for all  $SL(2, \CC)$-invariant vector fields $Y,Z$ tangents to $SL(2, \CC)/I$.
 
 Assume, by contradiction,  $S$ is locally modelled on $SL(2, \CC)/I$. Then $S$ inherits the   one form $\omega$ whose differential is a holomorphic
 volume form. This is in contradiction with Stokes's theorem. For further developpements of this method one can see~\cite{BL}
\end{proof}

{\bf Inoue surfaces}.
Recall that Inoue surface are compact complex surfaces in the class $VII_{0}$, which are not Hopf, and have  a vanishing second Betti number~\cite{Barth-Peters, Inoue}.

In~\cite{IKO, Klingler}  it is proved that any Inoue surface admits a (unique)   flat torsion-free holomorphic affine connection. Here we prove the following:

\begin{theorem}
Any holomorphic geometric structure $\tau$ on a Inoue surface is locally homogeneous.
\end{theorem}

\begin{remark} Inoue surfaces doesn't admit any non constant meromorphic functions~\cite{Barth-Peters}. 
\end{remark}

\begin{proof}

  Let  $\nabla_{0}$ be a  flat torsion-free   holomorphic  affine  connection on  $S$.  We prove that  the holomorphic rigid geometric structure $\tau'=(\tau, \nabla_{0})$ is locally homogeneous. 
  Since $a(S)=0$, theorem~\ref{result} implies $\tau'$ is locally homogeneous on some maximal open dense set $S \setminus E$, where $E$ is a compact analytic subset of $S$.

We want to show that $E$ is void.

We prove first that, up to a double cover of $S$, the subset $E$ is a smooth submanifold in $S$ (this is always true if $S$ is a finite number of points; but here $S$
might have  components  of complex dimension one).

Assume, by contradiction, that  $E$ is not a smooth  submanifold in $S$.

Choose $p \in E$ a singular point in  $E$. In particular, $p$ is not isolated in $E$, but $p$ is isolated among the singular points of $E$. Since $Is^{loc}(\tau')$ preserves $E$, it has to preserve the set of its singular points and thus it fixes $p$. Consequently  any local Killing field defined
in the neighborhood of $p$, has to vanish in $p$.

Denote by  $\mathcal G$ the Lie algebra of local Killing  fields in the neighborhood of $p$. Since  $\mathcal G$ acts transitively on an open set, its dimension is at least  $2$.

Each element of  $\mathcal G$ preserves  $\nabla$ and fixes 
 $p$. In exponential coordinates the $\mathcal G$-action in the neighborhood of $p$ is linear. This gives an embedding of  $\mathcal G$ in the Lie algebra of $GL(2, \CC)$. In particular,
$\mathcal G$ is of dimension  $\leq 4$.

Suppose first that $\mathcal G$  is of dimension $2$.  The corresponding isotropy subgroups of $GL(2, \CC)$ are conjugated either to the group of diagonal matrices, or to one of the following subgroups   $\left(  \begin{array}{cc}
                                                                a   &   b \\
                                                                 0     &  a^{-1} \\
                                                                 \end{array} \right)$, with  $a \in \CC^*$ and  $b \in \CC$, $\left(  \begin{array}{cc}
                                                                1   &   m \\
                                                                 0     &  n \\
                                                                 \end{array} \right)$, with  $m  \in \CC$ and  $n \in \CC^*$, or $\left(  \begin{array}{cc}
                                                                m'   &   n' \\
                                                                 0     &  1 \\
                                                                 \end{array} \right)$, with  $m' \in \CC^*$ and  $n \in \CC$.

    In the first case,  the invariant closed subset $E$  coincides, in exponential coordinates, with  the union of the two eigendirections. In the last two  cases, $E$ locally coincides with the invariant line $y=0$. In all  situations, up to a double  cover of $S$, the analytic set $E$ is smooth: a contradiction.

     We settle now the case where  $\mathcal G$ is of dimension  $3$ ou $4$. Then  the image of $\mathcal G$ by the isotropy representation in $p$  is conjugated  in $GL(2, \CC)$ to one of the following subgroups:
      $SL(2, \CC)$, $GL(2, \CC)$ or  the group of inversible upper triangular matrices. But $GL(2, \CC)$ and $SL(2, \CC)$ don't admit any invariant subset other 
     than $p$, which will be an isolated point in  $E$: impossible.
     
     In the last situation,  $E$ locally coincides, as before, with the unique invariant line and it is smooth.
     
     Up to a double cover, $E$ is a holomorphic submanifold in  $S$. If $E$ admits any component of dimension one, then  this component will be a union of closed curves. But Inoue surfaces doesn't contain any curve~\cite{Inoue}.
     
     This proves  $E$ is a finite number of points.  Assume, by contradiction,  that $E$ is not void and consider  $p \in E$.
     
     The previous arguments show that the Lie algebra  $\mathcal G$ is isomorphic either to  $sl(2, \CC)$, or to $gl(2, \CC)$. 
     
     Assume first  $\mathcal G =sl(2, \CC)$. As before, the local action of $\mathcal  G$ in the neighborhood of $p$ is conjugated to the action of  $\mathcal G$ on  $T_{p}S$ which coincides with the standard linear action of  $sl(2, \CC)$
     on  $\CC^2$.  This action has two orbits: the point $p$ and $\CC^2 \setminus \{p \} $.
     
    The stabilizer  $H$  in $SL(2, \CC)$ of a non zero vector  $x \in T_{p}S$  is conjugated to the following  one-parameter subgroup of 
     $SL(2, \CC)$:   $\left(  \begin{array}{cc}
                                                                1  &   b \\
                                                                 0     &  1 \\
                                                                 \end{array} \right)$, with  $b \in \CC$. Observe that the action of $G$ on $G/H$ preserves a non trivial holomorphic vector field. The expression of
   this vector field  in linear coordinates $(z_{1}, z_{2})$  in the neighborhood of $p$ is  $\displaystyle z_{1}     \frac{\partial}{\partial z_{1} }  +  z_{2}  \frac{\partial}{\partial z_{2}}$.

          Since  $S \setminus E$ is locally modelled on $(G,G/H)$,  this vector field  is well defined on $S \setminus E$. But, $E$ is of complex codimension $2$ in $S$ and the extension theorem of  Hartogs~\cite{GH}  shows that  the vector field  extends to a global non trivial holomorphic vector field $X$ on  $S$. The vector field $X$ is  $\mathcal G$-invariant on $S \setminus E$ and hence on   full $S$. Since the isotropy action of $SL(2, \CC)$ at $p$ doesn't preserve any non trivial vector in  $T_{p}S \simeq \CC^2$, it follows that $X(p)=0$. This is impossible since
         Inoue surfaces doesn't bear any {\it singular}  holomorphic vector field~\cite{Inoue}.
           
          The proof is the same in the case  $\mathcal G=gl(2, \CC)$.
      \end{proof}

      \subsection*{Meromorphic affine connections}
      
      In the following theorem we describe some  meromorphic affine connections on principal elliptic bundles over a compact Riemann surface of genus $g \geq 2$.

      \begin{theorem}  \label{meromorphic connections} Let $S$ be a principal elliptic bundle  over a compact Rieman surface $\Sigma$ of genus $g \geq 2$, with odd first Betti number.

     i)  A meromorphic  symmetric affine connection $\nabla$ on $S$ is invariant by the principal fibration if and only if the set of its poles intersects only a finite number of fibers. 
     
     ii) The space of those  previous connections admitting  simple poles on a single  fiber above some  point $\xi_{0} \in \Sigma$  is a complex affine   space of  dimension $5g+1$ and its underlying vector space identifies with  $H^0_{\xi_{0}}(K_{\Sigma}^2)  \times M_{2, \xi_{0}}^2$, where $H^0_{\xi_{0}}(K_{\Sigma}^2)$ is the vector space of meromorphic  quadratic differentials on $\Sigma$ with  a single pole of order at most two in $\xi_{0}$  and $M_{2,\xi_{0}}$ is the vector space of quasimodular forms of weight 2 on $\Sigma$ with a single simple pole at  $\xi_{0}$.
     
     iii) If  $\nabla$ is generic among the meromorphic affine connections which satisfy ii), then the Killing Lie algebra of $\nabla$ is generated by the fundamental
     vector field of the principal fibration.
     
       iv) The meromorphic affine connections which satisfy ii) are projectively flat away from the poles.
      \end{theorem}
      
      Recall that {\it a quasimodular form of weight $2$} on $\Sigma$ is a holomorphic map $f$ defined on the upper half-plane $H$ such that, for some $K \in \CC$, we have 
      $f(\xi)=f(\gamma \xi) (c_{\gamma} \xi +d_{\gamma})^{-2}  -K (c_{\gamma} \xi +d_{\gamma})^{-1}$,  for all $\gamma  =\left(  \begin{array}{cc}
                                                                a_{\gamma} &   b_{\gamma} \\
                                                                 c_{\gamma}     &  d_{\gamma} \\
                                                                 \end{array} \right)  \in SL(2, \RR)$ in the fundamental group of $\Sigma$~\cite{Azaiez, Zagier}. If $K=0$, then $f$ is a classical modular form. Theorem 9 in~\cite{Azaiez} shows that
  the space of  quasimodular forms of weight $2$ on $\Sigma$ which admits a simple pole in a single orbit is a complex vector  space of dimension $g+1$ (i.e. the quotient
  of the space of such quasimodular forms over those which are modular is one-dimensional).
      
      \begin{remark} The proof below shows that $S$ bears a (flat) holomorphic affine connection. Theorem~\ref{result} applies and, since any meromorphic function on $S$ is a pull-back of a meromorphic function on $\Sigma$~\cite{Barth-Peters},  the orbits of the pseudogroup of local isometries of any  meromorphic geometric structure on $S$ contain the fibers of $S$.
      \end{remark}
      
      \begin{proof} i) The projection of the set of poles of $\nabla$ on $\Sigma$ is a closed analytic set which is either a finite number of points, or full $\Sigma$. If this projection is $\Sigma$  and $\nabla$ is invariant by the principal fibration, then each point of $S$ is a pole, which is impossible. This proves the easy sense of the implication.
      
      We will describe now the space of connections $\nabla$,  such that   only a finite number of fibers contain poles of $\nabla$. Up to a finite unramified cover and a finite
      quotient, $S$ admits a holomorphic  affine structure (i.e. a flat symmetric  holomorphic affine connection $\nabla_{0}$) which can be built in the following way~\cite{Klingler}:
      
      Consider $\Gamma$ a discrete torsion-free subgroup  in $PSL(2, \RR)$ such that $\Sigma =\Gamma \backslash   H$, with $H$ the upper-half plane. Take any holomorphic projective structure on $\Sigma$, its developping map $\tau : H \to P^1(\CC)$ and its holonomy morphism $\rho : \Gamma \to PSL(2, \CC)$. This embedding of $\Gamma$ into $PSL(2,\CC)$ lifts to $SL(2, \CC)$ (this follows from the fact that orientable (real) closed $3$-manifolds have a trivial second Stiefel-Whitney class~\cite{MS}). Choose such a lift and consider $\Gamma$
      as a subgroup of $SL(2, \CC)$. Denote by $W = \CC^2 \setminus \{0 \}$ the $\CC^*$-tautological bundle over $P^1(\CC)$. The canonical affine structure of  $\CC^2$ induces a
      $\Gamma$-invariant affine structure on $W$ and hence a $\Gamma$-invariant holomorphic affine structure on  the pull-back $\tau^*(W) \simeq \CC^* \times H$. The $\Gamma$-action comes from the action by deck transformation on $H$ and from the $\rho$-action on $W$.
      
      The previous holomorphic affine structure on $\tau^*(W)$ is also invariant by homotheties in the fibers. Consider now $\Delta \simeq \ZZ$ a lattice in $\CC^*$ which acts by multiplication  in  the fibers of  $\tau^*(W)$
      and take the quotient of  $\tau^*(W)$ by $\Delta \times \Gamma$. The quotient is a principal elliptic bundle over $\Sigma$, with fiber $\Delta \backslash  \CC^*$, which identifies with $S$.
      
      The affine structure inherited by the universal cover $\CC \times H$ of $S$ is  the pullback of the previous affine structure of $\CC^* \times H$ by the map
      $$\CC \times H  \to \CC^* \times H$$$$(z, \xi) \to (e^z, \xi).$$
      
      The action of $\gamma=\left(  \begin{array}{cc}
                                                                a_{\gamma} &   b_{\gamma} \\
                                                                 c_{\gamma}     &  d_{\gamma} \\
                                                                 \end{array} \right)  \in \Gamma \subset  SL(2, \RR)$ on the universal cover $\CC \times H$ of $S$ is easily seen to be given by~\cite{Klingler}:
                                                                 
      $$\gamma (z, \xi)=(z+log(c_{\gamma}  \xi +d_{\gamma}), \gamma \xi), \forall (z, \xi) \in \CC \times H,$$ where $log$ is a determination of the logarithme and the $\gamma$ action
      on $H$ comes from the standard action of $SL(2, \RR)$ on $H$.
      
      The difference $\nabla - \nabla_{0}$ is a meromorphic $(2,1)$-tensor $\omega$ on $S$, or, equivalently, a $\Delta \times \Gamma$-invariant mermorphic (2,1)-tensor $\tilde{\omega}$ on the universal cover $\CC \times H$. Moreover, $\tilde{\omega}$ is invariant by $t(z, \xi)=(z +2i \pi, \xi)$.

Then we have,  $$\tilde {\omega}= f_{11}(z, \xi) dz \otimes dz \otimes \frac{\partial}{\partial z} + 
   f_{12}(z, \xi) dz \otimes d \xi  \otimes \frac{\partial}{\partial z}+
    f_{21}(z, \xi) dz \otimes dz \otimes \frac{\partial}{\partial \xi}+
     f_{22}(z, \xi) dz \otimes d \xi  \otimes \frac{\partial}{\partial \xi}+$$
    $$  g_{11}(z, \xi) d\xi  \otimes dz \otimes \frac{\partial}{\partial z}+
       g_{12}(z, \xi) d \xi  \otimes d \xi \otimes \frac{\partial}{\partial z}+
        g_{21}(z, \xi) d \xi  \otimes dz \otimes \frac{\partial}{\partial  \xi}+
         g_{22}(z, \xi) d \xi  \otimes d \xi  \otimes \frac{\partial}{\partial \xi},$$ with  $f_{ij}, g_{ij}$  meromorphic functions on $\CC \times H$, and $f_{ij}(\cdot, \xi), g_{ij}(\cdot, \xi)$ holomorphic
         except for a finite number of $\Gamma$-orbits in $H$.
         
         Notice that the difference between $\nabla_{0}$ and the standard affine structure of $\CC \times H$ is given by $f_{11}=f_{22}=g_{21}=1$, the others $f_{ij}, g_{ij}$ being trivial (see the straightforward computation in~\cite{Klingler}).
         
        Since  $\tilde \omega$ is $\Delta$-invariant and $t$-invariant, the functions $f_{ij}(\cdot, \xi)$ and $g_{ij}(\cdot, \xi)$ descend  on a   elliptic curve, for all $\xi \in H$. By the maximum principle, they
        are constant for all $ \xi \in H$ for which they are holomorphic. It follows that $f_{ij}(\cdot, \xi)$ and $g_{ij}(\cdot, \xi)$ are constant for $\xi$ lying in a  open dense subset of $H$ and, consequently,  for all $\xi \in H$. It follows that the functions $f_{ij}$ and $g_{ij}$ depend only on $\xi$ and 
        $\frac{\partial}{\partial z}$  is a Killing field for $\tilde{\omega}$. Since the fundamental generator of the principal fibration $\frac{\partial}{\partial z}$  preserves also $\nabla_{0}$, it is a Killing field for $\nabla$.

ii)         In the following we consider only symmetric connections  : $f_{12}=g_{11}$ and $f_{22}=g_{21}$.

                  The $\Gamma$-invariance of $\tilde{\omega}$  yields  to the following equations:\\

                 $(1)   f_{11}(\xi)=f_{11}(\gamma \xi) -c f_{21}(\gamma \xi)(c \xi +d)$
                  
                $ (2)     f_{12}(\xi)=f_{12}(\gamma \xi) (c \xi +d)^{-2}-2c^2f_{21}(\gamma \xi)  -c f_{22} (\gamma \xi) (c \xi +d)^{-1} +2c f_{11}(\gamma \xi) (c \xi +d)^{-1}$
                  
                  $(3)     f_{21}(\xi)=f_{21}(\gamma \xi) (c \xi +d)^2$
                  
                 $ (4)     f_{22}(\xi)=2f_{21}(\gamma \xi)c (c \xi +d)  +f_{22}(\gamma \xi)$
                  
                 $ (5)     g_{12}(\xi)= g_{12}(\gamma \xi) (c \xi +d)^{-4} +c^2 f_{11}(\gamma \xi)(c \xi +d)^{-2}   +  c  f_{12}(\gamma \xi) (c \xi +d)^{-3}-c^3f_{21}(\gamma \xi) (c \xi +d)^{-1}  -c^2f_{22}(\gamma \xi)(c \xi +d)^{-2}  -cg_{22}(\gamma \xi)(c \xi +d)^{-3}$
                  
                  $(6) g_{22}(\xi)=g_{22}(\gamma \xi)(c \xi +d)^{-2} +cf_{22}(\gamma \xi)(c \xi +d)^{-1} +c^2f_{21}(\gamma \xi),$

                   for all $\gamma =\left(  \begin{array}{cc}
                                                                a  &   b \\
                                                                 c   &  d\\
                                                                 \end{array} \right)$ in $\Gamma$.\\

           The equation $(3)$ implies that $f_{21}$ is a meromorphic  vector field on the compact Riemann surface $\Sigma$. Since a single pole of order at most one  is allowed, we get $f_{21}=0$ as a direct consequence of Riemann-Roch theorem and Serre duality~\cite{GH}.
           
           The equations $(1)$ and $(4)$ imply then that $f_{11}$ and $f_{22}$ are meromorphic  functions on $\Sigma$ with a single pole of order at most one. If $f_{11},f_{22}$
           are  not constant, they give a biholomorphism between $\Sigma$ and the projective line $P^1(\CC)$~\cite{GH}: a contradiction, since the genus of $\Sigma$ is $\geq 2$.
           
          We conclude that {\it  $f_{11}$ and $f_{22}$ are complex numbers.}
           
           Then we have :\\
           
           $(2')   f_{12}(\xi)=f_{12}(\gamma \xi) (c \xi +d)^{-2}  -c f_{22}  (c \xi +d)^{-1} +2c f_{11}(c \xi +d)^{-1}$
           
           $(5') g_{12}(\xi)= g_{12}(\gamma \xi) (c \xi +d)^{-4} +c^2( f_{11}-f_{22})(c \xi +d)^{-2}   +  c  (f_{12}-g_{22})(\gamma \xi) (c \xi +d)^{-3}$
           
           $(6')  g_{22}(\xi)=g_{22}(\gamma \xi)(c \xi +d)^{-2} +cf_{22}(c \xi +d)^{-1}$.\\
           
           It follows from $(2')$ and $(6')$  that $f_{12}$ and $g_{22}$ are quasimodular forms of weight two on $\Sigma$ with a single simple pole.  The space of those 
           quasimodular forms is a complex vector space of dimension $g+1$ (see~\cite{Azaiez}, theorem 9).

           Equation $(5')$ is equivalent  to the $\Gamma$-invariance of the quadratic
           differential $w(\xi)d\xi^2$, where
            $w=2g_{12}+ f'_{12} - g'_{22}$. It follows that $2g_{12}+ f'_{12} - g'_{22}$  is a meromorphic quadratic differential on $\Sigma$ with a single pole of order at most two.

           It is classicaly known (as an application of  Riemann-Roch theorem and Serre duality) that the space
           of quadratic differentials with a single pole of order at most $2$ is of complex dimension $3g-1$ (see, for example,~\cite{GH}).
           
   iii)  Let  $X=a(z, \xi) \frac{\partial}{\partial z}  + b(z, \xi) \frac{\partial}{\partial \xi}$ be a local holomorphic Killing field of a generic connection $\nabla$, in the neighborhood of a point where $\nabla$ is holomorphic( $a$ and $b$ are holomorphic local functions on $\CC \times H$).

  The equation of the Killing field is  $$\lbrack X, \nabla_{Y}Z\rbrack =\nabla_{\lbrack X,Y \rbrack}Z + \nabla_{Y}  \lbrack X, Z \rbrack,$$
for all  $Y,Z$ tangents to  $S$. It is enough to verify the equation for $(Y,Z)$ corresponding to  $(\frac{\partial}{\partial z}, \frac{\partial}{\partial \xi}),
(\frac{\partial}{\partial z}, \frac{\partial}{\partial z})$ et $(\frac{\partial}{\partial \xi}, \frac{\partial}{\partial \xi}).$

This yields to the following EDP system:\\

$ (1)~~a_{zz}  + (1+f_{11}) a_{z} + 2f_{12}b_{z}=0$

$(2)~~  b_{z z}  + (1+2f_{22}  -f_{11})b_{z}=0$

$\displaystyle (3)~~a_{z \xi}   + (f_{11}  -f_{22}) a_{\xi}  + g_{12}b_{z} +g_{11} b_{\xi}  + \frac{\partial f_{12}}{\partial  {\xi}}b=0$

$(4)~~b_{z \xi}  +(1+f_{22})a_{z} + (g_{22}-f_{12})b_{z}=0$

$\displaystyle (5)~~a_{\xi \xi} -g_{12}a_{z}  +(2f_{12}-g_{12})a_{\xi}  +2g_{12} b_{\xi}  + \frac{\partial g_{12}}{\partial \xi}b=0$

$\displaystyle (6)~~ b_{\xi \xi} + 2(1 +f_{22}) a_{\xi}  -g_{12}b_{z}+g_{22}b_{\xi} +  \frac{\partial g_{22}}{ \partial  \xi} b=0.$\\

The general solution of the first equation is  $b=\nu (\xi) e^{-\mu z}  + C(\xi)$, with  $\nu, C$ holomorphic functions of  $\xi$ and  $\mu = 1+2f_{22}  -f_{11}$ a complex constant.

We replace $b_{z}$ in the first equation and we get  $$\displaystyle a_{z}=\frac{\mu}{f_{11}-f_{22}} f_{12}(\xi) \nu(\xi)  e^{- \mu z} +A(\xi) e^{-(1+f_{11})z},$$ with $A$ a holomorphic function of  $\xi$.

Then  equation (4)  yields to 

$$\displaystyle \mu \lbrack - \nu'(\xi) + ( \frac{1+f_{11}}{f_{11}-f_{22}} f_{12} -g_{22}) \nu (\xi) \rbrack e^{-\mu z} + (1 +f_{22})A(\xi) e^{-(1+f_{11})z}=0.$$

For a generic $\nabla$, we have $f_{11} \neq f_{22}$, thus $\mu \neq 1+f_{11}$ and the functions $e^{- \mu z},e^{-(1+f_{11})z}$ are $\CC$-linearily independent. This 
implies $(1+f_{22})A(\xi)=0$ and, since for a generic  $\nabla$,  $f_{22} \neq -1$, we have  $A(\xi)=0$. 

We also  get  $$(I)~~ \nu'(\xi) = \nu (\xi) ( \frac{1+f_{11}}{f_{11}-f_{22}} f_{12} -g_{22}).$$

Now we check equation $(3)$. We replace the partial derivatives of $a$ and $b$ in $(3)$ and we get the following
$$\displaystyle \mu \lbrack \frac{f_{12}}{f_{11}-f_{22}} \nu'(\xi) - (g_{12} - \frac{f'_{12}}{f_{11}-f_{22}}) \nu(\xi)  \rbrack  e^{-\mu z} +(f_{12}C)'=0.$$

Since generically $\mu \neq 0$, the functions $e^{- \mu z}$ and $1$ are linearily independent, which yields to
$(f_{12}C)'=0$ and to

$$\displaystyle (II)~~\nu'(\xi) = \nu(\xi) \frac{1}{f_{12}} \lbrack (f_{11}-f_{22})g_{12}   -f'_{12}   \rbrack .$$

Relations  $(I)$ et $(II)$ are  compatible, for a generic connection, only if $\nu=0$. This implies $b=C(\xi)$, $a=B(\xi)$.

Our EDP system becomes:\\

$(3')~~(f_{11}-g_{21})a'  +f'_{12}b+f_{12}b'=0$

$(5')~a'' +(2f_{12}-g_{12})a'+   2g_{12}b'-g'_{12}b=0$

$(6')~~b''  +2(1+2g_{21})a'+g_{22}b'+g'_{22}b=0.$\\

Since $(f_{12}b)'=0$, we have $a'=0$ and thus $a$ is a constant function. Equation $(5')$ implies then $2g_{12}b'-g'_{12}b=0$, which, for a generic $\nabla$,  is compatible 
with $f_{12}b'+f'_{12}b=0$ only if $b=0$.

It follows that $X$ is a constant multiple of $\frac{\partial}{\partial z}$.

iv) The projective connection associated to the affine connection $\nabla$ is given by the following  second order ODE~\cite{Cartan}:

$$\xi''=K^0(z, \xi) + K^1(z, \xi) \xi' + K^2(z, \xi)(\xi')^2+ K^3(z, \xi)(\xi')^3,$$ where
$K^0=- f_{21}=0$, $K^1=(1+f_{11})-2(1+ f_{22})$, $K^2=-(g_{22}-2f_{12})$ et $K^3=g_{12}$.

Liouville~\cite{Liouville}, followed by Tresse~\cite{Tresse} and Cartan~\cite{Cartan},  proved that this projective connection is projectively flat if and only if both following invariants vanish:

$$L_{1}=2K^1_{z\xi}-K^2_{zz}-3K^0_{\xi \xi}-6K^0K^3_{z} -3K^3K^0_{z} + 3K^0K^2_{\xi}+3K^2K^0_{\xi}+
K^1K^2_{z}-2K^1K^1_{\xi},$$ $$L_{2}=2K^2_{z \xi}  -K^1_{\xi \xi}  -3K^3_{z z} + 6 K^3  K^0_{\xi}  + 3 K^{0}K^3_{\xi}-3K^3K^1_{z}  -3K^1K^3_{z}   -K^2K^1_{\xi} +2K^2K^2_{z}.$$

Here  $K^0=0$, $K^1$ is a constant function   and  $K^2, K^3$ depend only of  $\xi$. This implies the vanishing of both invariants $L_{1}$ and $L_{2}$.
 \end{proof}
        
 {\it I would like to thank G. Chenevrier and  G. Dloussky  for helpful conversations.}

${}$ 
\end{document}